\documentclass[11pt]{amsart}
\usepackage[dvips]{graphicx,color}   
\usepackage{amsmath}
\usepackage{amscd}
\usepackage{amsfonts,latexsym,amssymb}

\newcommand{\s}{\sigma}

\renewcommand{\o}{\omega}

\DeclareFontFamily{OT1}{rsfs}{}
\DeclareFontShape{OT1}{rsfs}{n}{it}{<->rsfs10}{}
\DeclareMathAlphabet{\curly}{OT1}{rsfs}{n}{it}

\newcommand{\C}{\mathbb{C}}

\newcommand{\N}{\mathbb{N}}
\newcommand{\R}{\mathbb{R}}

\newcommand{\CP}{\mathbb{CP}}

\newcommand{\sref}{s_{k,x}^{\text{\scriptsize ref}}}

\newcommand{\x}{\times}
\newcommand{\ox}{\otimes}

\newtheorem{proposition}{Proposition}[section]
\newtheorem{theorem}[proposition]{Theorem}
\newtheorem{definition}[proposition]{Definition}
\newtheorem{lemma}[proposition]{Lemma}

\newtheorem{corollary}[proposition]{Corollary}
\newtheorem{remark}[proposition]{Remark}

\def\>#1{{\bf #1}}             

\def\Gr{\hbox{{\rm Gr}}}
\def\Ker{\hbox{{\rm Ker}}}

\begin{document}

{\mbox{   }}

\vskip 2cm

\centerline{\Large\bf SUBMANIFOLDS OF SYMPLECTIC MANIFOLDS}

\bigskip

\centerline{\Large\bf WITH CONTACT BORDER}

\medskip

\vskip 2cm

\centerline{F. Presas}

\vskip 1.5cm

\medskip
\centerline{ {\it Departamento de \'Algebra, Universidad 
Complutense de Madrid, 28040 Madrid, Spain.}  }

\vskip 2cm

\centerline{\sc Abstract}
{\small We construct symplectic submanifolds of symplectic manifolds
with contact border. The boundary of such submanifolds is shown to be a
contact submanifold of the contact border.
We also give a topological characterization of the constructed
submanifolds by means of a ``relative Lefschetz hyperplane Theorem''.
We sketch some of the applications of the results.}

\vskip 1.5cm

\newpage

\section{Introduction.} \label{introduction}
A number of works have been developped from the foundational paper \cite{Do96}
that exploit the idea of ampleness in the symplectic and contact category.
The key idea has been to adapt the concept
of ``linear system'' to these cases. The techniques have 
provided a new insigth in symplectic topology, giving as byproduct
new symplectic invariants \cite{Au99b, Do99}. In the contact case, it
has been possible to mimic the symplectic counterpart \cite{IMP99}. The results contained
in \cite{Pr00} open the way for constructing contact invariants
although it will be less direct because, contrary to the symplectic situation,
we do not have canonical constructions (up to symplectic isotopy). 

The aim of this paper is to put together the two constructions, in the symplectic
and contact category to push-forward the progress in the study
of symplectic submanifolds with contact border. This concept appears
naturally when defining ``cobordisms'' in the symplectic category. 

We call sym-con category the category defined by symplectic manifolds with contact
border. An object in this category is a set $(M,\omega, C,\theta)$ such that $(M,\omega)$
is an open symplectic manifold with compactification $\bar{M}=M\cup C$ and 
$(C,\theta)$ is a cooriented
contact manifold whose structure is compatible with the symplectic one in the usual sense
(see Subsection \ref{basic} for details). Along the paper the dimension on the symplectic
manifold $M$ will be $2n+2$, except when an explicit mention is made. The hermitian 
line bundle $L$ whose curvature is $-i\omega$ will be called prequantizable line bundle.

In this article we study the possibility of constructing submanifolds in this
sym-con category with topological properties similar to divisors in complex
projective geometry. We provide a complete topological characterization of
these submanifolds which are approximately holomorphic in the sense of
\cite{Do96,Au97,IMP99}. The main result of this paper is
\begin{theorem} \label{main_thm}
Let $(M,\omega)$ be a symplectic manifold of integer class with prequantizable
bundle $L$ and with contact border
$(C,\theta)$. Fix a rank $r$ complex vector bundle $E$ over $M$.
For $k$ large enough, there exists a symplectic submanifold $W$ of $M$,
which is Poincar\'e dual of $c_r(L^{\ox k}\otimes E)$,
satisfying that $\bar{W}\cap C$ is a contact submanifold of $C$. Moreover
the inclusion $i: \bar{W}\to \bar{M}$, induces an isomorphism in relative homology groups 
through the natural morphism $H_p(\bar{W}, \bar{W}\cap C)\to H_p(\bar{M},C)$ for 
$p<n-r$ and an epimorphism for $p=n-r$.
\end{theorem}
First we do notice that Theorem \ref{main_thm} could follow by a more or less
straightforward combination of ideas in \cite{Au97} and \cite{IMP99}. 
For this one might only to extend the contact submanifolds
constructed in \cite{IMP99} to a small neighborhood of $C$ in $M$. However, this kind
of approach presents problems of difficult solution. We detail a little more this
question, using the approximately holomorphic tools, in Subsection \ref{asint_holo}
( cfr. Remark \ref{triviality}).

So we have chosen an alternative way to attack the problem. We will define directly 
global sections which solve the problem near the border. For this we need
to revise the local theory developed in \cite{IMP99, Pr00}. The solution goes
through the use of approximation theory and a refined Jackson's theorem.
Jackson's theorems \cite{Ch66} provide boundings for the error
made when we approximate a differentiable function by a polynomial of a given
degree in terms of the derivatives of the function. This kind of results were
used by S. Donaldson implicitly in the foundational work \cite{Do96},
when he approximated asymptotically holomorphic functions by polynomials.
But now, we need a similar result for any function such that we only control the
norm of the derivatives.

In Section \ref{defini} we will state the approximately holomorphic theory in the 
sym-con category, following the notations of \cite{Do99, Pr00}. We reduce the proof
of Theorem \ref{main_thm} to a transversality result in the border. In Section 
\ref{topol} we prove
the relative Lefschetz hyperplane theorem as stated in Theorem \ref{main_thm} and, also
we characterize the Chern classes of the constructed submanifolds. 

We must stress that the topological results are new even in the integrable complex case 
and offer a new insight in the topological structure of the sym-con manifolds, showing
again that the idea of Eliashberg of studying contact manifolds as the more natural
border definition in the symplectic and k\"ahler category is powerful and
rich in consequences.

{\bf Acknowledments:} \newline
I am very grateful to S. Donaldson and D. Auroux by his kindness passing me their 
preprints \cite{Do99} and \cite{Au99b}. 
I want to thank to the members of the GESTA \footnote{Geometr\'{\i}a simpl\'ectica con
t\'ecnicas algebraicas, CSIC-UC3M, 2000.} seminar in Madrid their support and 
interest through the
elaboration of this work. I want to thank especially to Vicente Mu\~noz by his
careful reading of this paper and his useful suggestions about the topological results.

\section{Definitions and results.} \label{defini}
We state along this Section the basic notions in symplectic and contact topology
needed in what follows. Also we state the main result in terms of
the tools introduced and sketch the idea of the proof of Theorem \ref{main_thm}.

\subsection{Sym-con manifolds.} \label{basic}

The symplectization $S_D(C)$ of a contact manifold $(C,D)$ is defined as
$$S_D(C)= \{ \theta\in T^*C: \Ker ~ \theta= D(\pi (\theta)) \},$$
where $\pi: T^*C \to C$ is the standard projection. It is easy to check that this
manifold has a canonical exact symplectic structure provided by the exterior
differential of the Liouville $1$-form
$$ \alpha(v_{\theta})= \theta(\pi_*(v)), \forall v_{\theta}\in TS_D(C). $$
The symplectization $S_D(C)$ has structure of a $\R^*$-principal bundle over $C$. We are 
particularly interested on the cooriented case (called exact case as well). The contact
manifold is said to be cooriented if there exists a global $1$-form $\theta$ in $C$
satisfying that $\Ker ~ \theta= D$. In that case the form $\theta$ provides a section
of $S_D(C)$ and this becomes a trivial bundle. Fixing a $1$-form $\theta$
we can identify canonically $S_D(C)=C\times (-\infty,0)\cup (0, \infty)$. 
In this paper we will only use the cooriented case and from now on
we call symplectization to the connected component $C\times (0, \infty)$
instead of the total set. We can give an explicit formula for the symplectic
form in that manifold, we set up the following (canonical) isomorphism:
\begin{eqnarray*}
S_D(C) & \to & C\times (0, \infty) \\
\lambda\cdot \theta(x) & \to & (\pi(\theta(x)), \lambda).
\end{eqnarray*}
Then we obtain
$$d\alpha= d\lambda \wedge \pi^*\theta +\lambda\cdot \pi^*d\theta. $$
The contact manifold $C$ can be embedded in the symplectization
through the graph of $\theta$, namely as the contact hypersurface 
$\hat{C}_{\epsilon}=\{ (p,\lambda)\in C\x \R^+: \lambda=\epsilon \}$.
This embedding will be called the $\epsilon$-embedding of $C$ in $S_D(C)$,
denoted as $i_{\epsilon}$. If there is not risk of confussion, we usually denote $\hat{C}_1$
by $\hat{C}$. Recall that through this family of embeddings the distribution $D$ defines
a distribution $\hat{D}$ of $2$-codimensional spaces in $S_D(C)$, it is obvious that
$\hat{D}$ is symplectic with respect to the canonical symplectic structure
in the symplectization.

Recall that a symplectic manifold $(M,\omega)$ has a contact border $C$ when
we can identify through a symplectomorphism a neighborhood $V$ of the border $C$ with one of 
the two standard models: $C\times [a,b)$ or $C\times (a,b]$, for some $a,b\in \R^+$, with
$a<b$. In the first case we will say that the manifold has concave border and in the second
one convex border. Moreover we can generalize the definition to include mixed
cases. So in general a symplectic manifold has contact borders $(C_1^1,C_1^2,
\ldots, C_1^{cc}, C_2^1, \ldots, C_2^{cv})$, where each $C_i^j$ is a connected
contact cooriented manifold, if we can decompose a neighborhood $V$ of the border
into $cc+cv$ connected components $V_1^i$, $V_2^j$ ($i=1, \ldots, cc$ and $j=1, \ldots cv$)
such that each of the $V_1^i$ is symplectomorphic to the local concave model
defined by $C_1^i$ and respectively with $V_2^j$ and the convex model.

An important observation is that convex and concave models are not
equivalent and provide very different problems in the sym-con category
(see i.e. the pseudo-holomorphic curves construction in \cite{El98}).

\begin{remark}
In the literature, definitions above are called sometimes ``strictly convex 
(or concave) contact borders'' to distinguish them from a weaker notion defined
by Eliashberg as follows. We will say that a symplectic manifold $(M, \omega)$ has
``weak'' contact borders  $(C_1^1,C_1^2,
\ldots, C_1^{cc}, C_2^1, \ldots, C_2^{cv})$, where each $C_i^j$ is a connected
contact cooriented manifold, if we can decompose a neighborhood $V$ of the border
into $cc+cv$ connected components $V_1^i$, $V_2^j$ ($i=1, \ldots, cc$ and $j=1, \ldots cv$)
such that each of the $V_1^i$ is diffeomorphic to the local concave model
defined by $C_1^i$ (respectively with $V_2^j$ and the convex model). The 
diffeomorphism $\phi: V \to C\times [a,b)$ must satisfy that $i_{a}^*\phi_*\omega$
is non degenerated when restricted to the distribution $D$. (resp.
in the convex model). All the theory developped in this article can be adapted
with slight modifications to this more general case, obtaining submanifolds
with ``weak contact border''. We do not detail this along the article but
the reader can translate all the proofs to that case.
\end{remark}

As in the closed manifold case we can construct a complex line bundle, $L$,
over a symplectic manifold $(M,\omega)$ whose curvature form is $-i\omega$,
provided an integrality condition is satisfied, namely $[\omega/2\pi]$ has to be
the lifting of an integer class. This
bundle is usually called a prequantizable bundle, because of the geometric
quantization setting. Moreover the precedent considerations assure
that the bundle extends to the border defining a line bundle whose
curvature form is $-id\alpha$ (under the standard models identifications) in each 
connected component. We will denote by $L$ the prequantizable bundle in $M$ and 
also its extension to $C$.

\begin{definition}
A sym-con manifold $(M,\omega, C,\theta)$ is a symplectic manifold $(M,\omega)$
with compactification $\bar{M}$ satisfying that $\bar{M}-M=C$ admits a contact
structure $\theta$ which defines a contact border for $M$.
\end{definition}

The following trivial result is the symplectic analogue of the connected
sum theorem in topology. 
\begin{lemma} \label{tubo}
Given two sym-con manifolds $(M_1, \omega_1, C, \theta)$
and $(M_2, \omega_2, C, \theta)$ with connected convex and concave borders respectively.
Then, for a suitable nonzero constant $\lambda$ ,
there exists a closed symplectic manifold $(\tilde{M}, \omega)$ and two
symplectic embeddings $\varphi_1:(M_1, \omega_1) \to (M, \omega)$ and
$\varphi_2:(M_2, \lambda \omega_2) \to (M, \omega)$ satisfying that
$\varphi_1(M_1)\cup \varphi_2(M_2)=M$.
\end{lemma}
The manifold $M$ is usually denoted as $M_1 \cup_C M_2$ and topologically is a
connected sum along $C$.

\noindent {\bf Proof:}
We have only to use the standard models of the borders to glue the manifolds
symplectically. Say that near the border the local model for 
$M_1$ is $C\times(a_1,b_1]$ and for $M_2$ is $C\times[a_2,b_2)$. If we find
that $(a_1,b_1]\bigcap [a_2,b_2)\neq \emptyset$ then we are finished.
If not we substitute the symplectic form $\omega_2$ by $\lambda\omega_2$.
This produces a change in the local model of $M_2$ which is now $C\times[\lambda a_2,
\lambda b_2)$. Obviously a suitable choice of $\lambda$ reduces the problem to
the precedent one.
\hfill $\Box$

Recall from the proof that it is not very important that the contact
forms chosen in the two borders coincide, if the distribution is the same. 
In fact, the symplectic connected sum along
a contact border does not depend on this choice, because the symplectic structure
of the symplectization does not depend on the choice of contact form.

We can always add a symplectic collar in the border of a sym-con manifold. This is
the content of the following
\begin{corollary} \label{collar}
Let $(M, \omega, C, \theta)$ a sym-con manifold, then we can find a manifold
$(M', \omega', C, \theta)$ such that there exists a symplectic embedding of $(M, \omega)$
in $(M', \omega')$ satisfying that the compactification of $M$ does not intersect
the border $C$ of $M'$.
\end{corollary}
\noindent {\bf Proof:}
It is a direct application of the precedent Proposition choosing $M_1=M$
and $M_2=C\times[1/2,3/2)$ if the border of $M$ is convex (resp. $M_2=
C\times(1/2,3/2]$ if the border is concave).
\hfill $\Box$

\begin{definition}
A contact hypersurface $(C,\theta)$ in a symplectic manifold $(M,\omega)$ is a
hypersurface in $M$ supporting a $1$-form $\theta$ such that $D=\Ker ~ \theta$
is a contact distribution and $d\theta_{|D}=\omega_{|D}$.
\end{definition}
If, using Corollary \ref{collar}, we add a symplectic collar to a sym-con
manifold $(M,\omega, C, \theta)$ then, the submanifold
$C$ is a contact hypersurface in the enlarged manifold $M'$. We will use this 
idea afterwards.

We define a compatible almost-complex structure $J$ in a sym-con manifold $(M,\omega,
C,\theta)$ as a compatible almost-complex structure in $(M,\omega)$ such that the restriction
of $J$ to the contact border $C$ leaves invariant the distribution $D$. By using
the local model it is obvious that in this case the restriction of $J$ to the distribution
$D$ provides a compatible almost-complex structure in the symplectic bundle $D$.
It is easy to check that the moduli space of such structures is contractible. For this 
we use the same arguments that in the symplectic and in the contact case.

As always, when we fix a compatible almost-complex structure, we automatically
obtain a metric $g$ on the manifold $(M, \omega)$ as $g(v,w)=\omega(v,Jw)$.
We refer to this metric as the symplectic metric. We define also the
$k$-rescaled symplectic metric as $g_k=kg$.

\subsection{Contact manifolds.}
Now, we recall some basic ideas about contact geometry. We assume that
$(C,D)$ is a cooriented contact manifold where we have fixed a contact form $\theta$.
This contact form determines a vector field $R$ by the conditions:
$$ i_R \theta=1, ~~ i_R d\theta= 0, $$
which is called the Reeb vector field. As in the symplectic case when we fix a compatible
almost-complex structure $J$ we obtain a metric in the contact manifold as
$g(v,w)=\theta(v)\theta(w)+d\theta(v,Jw)$, which is called the contact metric. The
$k$-rescaled contact metric is defined as $g_k=kg$. We are abusing notation by using the
same letter to denote the symplectic and contact metrics, but it is easy to check that
in a sym-con manifold the restriction of the $k$-rescaled symplectic metric to the contact 
border coincides with the precedent definition, justifying our notation. However, an 
important
change of behaviour appears in the contact case. Namely the $k$-rescaled symplectic
metric is the symplectic metric associated to the form $k\omega$, but in the
contact case the $k$-rescaled contact metric is not the contact metric associated
to $k\theta$. This difference is fundamental to develop the theory and will
reflect, in the contact case, the localization process which appears in Donaldson's theory.
We formalize this idea with the following definitions.

\begin{definition}
The maximum angle between two subspaces $U,V\in \Gr_{\R}(r,n)$ is defined as
$$ \angle_M(U,V)=\max_{u\in U} \max_{v\in V} \angle(u,v). $$
\end{definition}
This angle defines a distance in the topological space $\Gr_{\R}(r,n)$ (for details see
\cite{MPS99}).

\begin{definition}
Let $D_k$ be a sequence of contact distributions in $\R^{2n+1}$. The sequence is
called $c$-asymptotically flat in the set $U\in \R^{2n+1}$ if
$$ \angle_M(D_k(0),D_k(x))\leq ck^{-1/2}, \ \mbox{\rm for all} \ x\in U. $$
The sequence is called asymptotically flat if there exist some $c$ for which
it is $c$-asymptotically flat.
\end{definition}

The standard contact structure in $\R^{2n+1}$ is defined as $\theta_0= ds+\sum_{j=1}^n 
x_jdy_j$, where $(x_j,y_j,s)\in \R^{2n+1}$. The $k$-rescaled contact metric is the
contact metric associated to $\theta_{k^{1/2}}= k^{1/2}ds+\sum_{j=1}^n x_jdy_j$,
which is obtained from $\theta_0$ scaling the coordinates by a factor $k^{1/2}$.
So, it is clear that, at any small neighborhood of a given point, when we apply
the set of metrics $g_k$ we obtain as a result, passing to a fixed Darboux trivialization,
a sequence of contact forms $\theta_{k^{1/2}}$, by scaling with a factor $k^{1/2}$
in $\R^{2n+1}$, which is obviously asymptotically
flat in any bounded set in $\R^{2n+1}$.

Given any asymptotically flat sequence of distributions $D_k$ in $\R^{2n+1}=\C^n\times \R$,
satisfying that $D_k(0)=\C^n \times \{0 \}$,
there exists a canonical almost-complex structure in a neigborhood of the origin, for $k$
large enough. We only have to lift the complex structure defined in $\C^n$ to the 
distribution $D_k$ using the pull-back of the vertical projection (which is an isomorphism
near the origin for $k$ large enough).

Finally if we have a contact hypersurface $(C,\theta)$ in a symplectic manifold 
$(M,\omega)$ we can choose
a compatible almost-complex structure which makes the distribution $D$ $J$-invariant.
In fact, in this case we can identify symplectically a neighborhood of $C$ with
a neighborhood of the $1$-embedding of $C$ in the symplectization $S_D(C)$ and the 
almost-complex structure can be chosen to make the distribution $\hat{D}$ $J$-invariant 
in this neighborhood (through the identification). This kind of almost complex
structures will be called compatible with the hypersurface.
Suppose that we have added a symplectic collar to a sym-con manifold. A compatible 
almost-complex structure $J$ in the sym-con manifold admits an extension to an almost
complex structure $\tilde{J}$ in the enlarged manifold which is compatible
with the contact hypersurface $C$.

\subsection{Sequences of bounded sections.} \label{asint_holo}
We recall from \cite{Do99, Pr00} the approximately holomorphic setting.
We adapt it to our present work and ideas. A uniform constant, polynomial, etc.
is a constant, polynomial, etc which does not depend on the chosen point
of the sym-con manifold nor the integer $k$ appearing in the context.

Now, we introduce the notion of asymptotically holomorphic sections which is
one of the key points. All the norms in the definitions to follow
are defined with respect to the sequence of metrics $g_k$.

\begin{definition}[\cite{Do99}]
A sequence of sections $s_k$ of the hermitian bundles $E_k$ over the symplectic
manifold $(M,\omega)$ has $C^r$-bounding $c$ at the point $x$ if it satisfies
\begin{eqnarray*}
|s_k(x)|& <& c, \\
|\nabla^j s_k(x)| & < & c, ~~ \forall j=1, \ldots, r, \\
|\nabla^{j-1} \bar{\partial} s_k(x)| & < & ck^{1/2}, ~~ \forall j=1, \ldots, r.
\end{eqnarray*}
The sequence has uniform $C^r$-bounding $c$ if it satisfies these boundings at every point.
\end{definition}

\begin{definition}[\cite{Pr00}]
A sequence of sections $s_k$ of the hermitian bundles $E_k$ over the contact
manifold $(C,\theta)$ has mixed $C^r$-boundings $(c_D,c_R)$ at the point $x$ if it satisfies
\begin{eqnarray*}
|s_k(x)|& <& c_D, \\
|\nabla^j_D s_k(x)| & < & c_D, \forall j=1, \ldots, r, \\
|\nabla^j_R s_k(x)| & < & c_R, \forall j=1, \ldots, r, \\
|\nabla^{j-1} \bar{\partial} s_k(x)| & < & c_Rk^{-1/2}, \forall j=1, \ldots, r.
\end{eqnarray*}
The sequence has uniform mixed $C^r$-boundings $(c_D,c_R)$ if it satisfies 
these boundings at every point.
\end{definition}

\begin{definition} \label{global_b}
A sequence of sections $s_k$ of the hermitian bundles $E_k$ over the sym-con
manifold $(M,\o,C,\theta)$ has global $C^r$-boundings $(c,c_D,c_R)$ at a point $x\in C$ if 
$s_k$ restricted to $M$ has $C^r$-bounding $c$ and restricted to
$C$ has uniform mixed $C^r$-boundings $(c_D,c_R)$. 
\end{definition}
As usual the $C^r$-openness is important in this kind of definitions. Namely,
if we have sections $s_k^1$ and $s_k^2$ with global $C^r$-boundings $(c^1,c_D^1,c_R^1)$
and $(c^2,c_D^2, c_R^2)$ then $s_k^1+s_k^2$ has global $C^r$-bounds
$(c^1+c^2, c_D^1+c_D^2, c_R^1+c_R^2)$. An analogous property is satisfied by the
other types of boundings.

Definition \ref{global_b} also applies to sequences of sections defined over a
symplectic manifold $M$ which contains a contact hypersurface $C$, being the 
definition in this case the natural one.

The other key ingredient is the notion of transversality with estimates. We set up
it in a general way following \cite{IMP99}.

\begin{definition} \label{transD} 
Let $s$ be a section of a complex vector bundle $E$ over the Riemannian
manifold $X$ with distribution $D$, and
$\eta > 0$.  The section $s$ is said to be $\eta$-transverse to {\bf 0}
along $D$ at a point $x\in X$ if it is satisfied at least one the following conditions
\begin{enumerate}
\item $|s (x)| < \eta$,
\item the covariant derivative restricted to $D$,
$\nabla_D s \colon D_x\subset T_xX \to E_x$, is surjective and has a right 
inverse of norm less than $\eta^{-1}$.
\end{enumerate}
The section is $\eta$-transverse to {\bf 0} along $D$ in a set $U$, if it is $\eta$-transverse
at all the points of $U$.
\end{definition}
In the symplectic case $D=TM$ and in the contact case $D$ is the contact distribution.
This definition is $C^1$-open in the sense that there exists a constant $c_0$, only depending
in the dimensions, such that if $s$ is $\epsilon$-transverse to {\bf 0} along $D$ and
$|\s-s|<\alpha$ then $\s$ is $(\epsilon-c_0\alpha)$-transverse to {\bf 0} along $D$.
It is possible to precise a little more in the contact case. Namely, 
again there exists a constant $c_0'$ such that if $s$ is $\epsilon$-transverse to {\bf 0} 
along $D$ and $|\s-s|$ has mixed $C^1$-boundings $(\alpha, c_R)$ then $\s$ is 
$(\epsilon-c_0'\alpha)$-transverse to {\bf 0} along $D$.

With these definitions at hand we reduce the proof of Theorem \ref{main_thm}
to the following:
\begin{proposition} \label{main_propo}
Let $(M,\omega)$ be a symplectic manifold of integer class with contact border
$(C,D)$. Let $U$ a compact set in $M$ which does not intersect $C$. Fix a rank $r$ 
complex vector bundle $E$ over $M$.
Fix a constant $\epsilon>0$ and a compatible almost-complex structure in the sym-con
manifold. Let $s_k$ be a sequence of sections with global $C^3$-boundings of the
bundles $E\ox L^{\ox k}$. Then there exists a sequence of sections $\s_k$ with
global $C^3$-boundings such that $|s_k-\s_k|_{C^1,U}<\epsilon$ and satisfying that
$\s_k$ is $\eta$-transverse to {\bf 0} in $M$ and $\eta$-transverse to {\bf 0}
along $D$ in $C$.
\end{proposition}
Observe that near the border we cannot control the $C^0$-norm of the perturbation. \newline
{\bf Proof of the existence part of Theorem \ref{main_thm}:}
Take a sequence of sections $\sigma_k$ given by Proposition \ref{main_propo}.
We only have to apply Lemma 5 of \cite{IMP99} to the manifold
$M$ and to the border $C$ respectively to obtain that the zero sets
are symplectic and contact. Only notice that the asymptotically holomorphic
sequences of that article correspond to our $C^r$ and mixed $C^r$-boundings.
\hfill $\Box$

\begin{remark} \label{triviality}
A direct approach for proving Proposition \ref{main_propo}
should be to define a mixed $C^r$-bounded sequence of sections in the border 
$C$ which be $\eta$-transverse to {\bf 0} provided by \cite{IMP99,Pr00} 
and try to extend the sequence to the symplectic manifold. In
fact, the boundings in the derivatives work to produce this extension and the 
holomorphicity condition gives us the 
derivatives of a given section in the normal direction. But making the computations
in detail we find that we are able to extend the construction in an asymptotically
holomorphic way only to a strip of $g_1$-radius $O(k^{-1/2})$ from the border. This is not
enough to multiply by a cut-off function and so to define the section all over $M$, because
the global boundings are destroyed. We would need a strip
of $g_1$-radius $k^{-1/3}$, but the arrangement to get the boundings in this strip is 
not clear. In the next paragraphs we explain the method of proof that we use to overcome
this difficulty.
\end{remark}

\subsection{Proof of Proposition \ref{main_propo}.}
We state from the results of \cite{Au97} the following 
\begin{theorem}[Adaptation of Theorem 2 in \cite{Au97}] \label{Auroux}
Let $E$ be a complex vector bundle of rank $r$ over a symplectic manifold $(M,\omega)$
of integer class (not necesarilly compact). Let $J$ be a compatible 
almost-complex structure. Fix a constant $\epsilon>0$ and a compact set $U$ in $M$, 
and let $s_k$ a sequences of sections with $C^r$-bounding $c$ of the bundles $E\ox L^{\ox k}$.

Then there exists a uniform constant $\eta>0$ (depending only on $\epsilon$ and $c$)
and a sequence $\s_k$ of sections with $C^r$-bounding $\epsilon$ such that $s_k+\s_k$
is $\eta$-trasnverse to {\bf 0} over $U$.
\end{theorem}
{\bf Proof:}
The only difference with respect to Auroux' result is that we do not impose
the closedness of the manifold $M$. But Auroux techniques are purely local. So there is no
reason to impose the closedness of the manifold $M$. The only important point is to 
guarantee the compactness condition and this is assured by restricting ourselves to a 
compact set $U\subset M$.
\hfill $\Box$

The existence of the border makes impossible to set up the $1$-parametric discussion
of \cite{Au97} as is shown in \cite{IMP99, Pr00}.
We want to reduce the proof of Proposition \ref{main_propo} to the following result
\begin{theorem} \label{trans_bor}
Let $\epsilon>0$, $\alpha>0$. 
Given a cooriented contact manifold $(C,\theta)$ and given the 1-embedding
of $C$ in the symplectization $S_D(C)=C\times \R$, fix a complex vector bundle $E$ over
the symplectization. Then given a global $(c,c_D,c_R)$-bounded in $C^r$-norm sequence 
of sections $s_k$ of the bundles $E\otimes L^{\otimes k}$, there exists another
sequence of sections $\tau_k$ with global $(c', \epsilon, c_R')$ $C^r$-boundings 
satisfying, for $k$ large enough, that
\begin{enumerate}
\item $\tau_k$ is supported in $C\times (1-\alpha, 1+\alpha)$.
\item The restriction to $\hat{C}$ of $s_k+\tau_k$ is $\eta$-transverse to {\bf 0}
along the distribution $D$ in $\hat{C}$, for some uniform constant $\eta>0$.
\end{enumerate}
\end{theorem}
We assume this result, which will be proved in Section \ref{borders} and then we prove:

\noindent {\bf Proof of Proposition \ref{main_propo}:}
Fix a compatible almost-complex structure in the sym-con manifold.
Enlarge, adding a symplectic collar, the sym-con manifold $(M,\omega,C,\theta)$ to 
obtain a new symplectic manifold $M'$ where $C$ is a closed contact hypersurface.
For $\alpha>0$ small enough we can identify symplectically a neighborhood $V$ of $C$ with 
the neighborhood $C\times (1-\alpha, 1+\alpha)$. We can extend the almost-complex
structure with one compatible with the hypersurface. Fix a sequence of sections $s_k$
with uniform $C^r$-bounding $c$ in $(M',\omega')$, obviously $s_k$ has global $C^r$-bounds
$(c,c,c)$ in the initial manifold $M$. 

Then we apply Theorem \ref{trans_bor}
to $V\simeq C\times(1-\alpha,1+\alpha)$ perturbing the sequence $s_k$ to obtain a
new sequence $\s_k$ which is $\eta$-transverse along the distribution $D$ on
$C$. To finish we perform a perturbation $\tau_k$ of $C^1$-norm less than $\frac{\eta}
{2c_o}$, where $c_o$ is the constant of $C^1$-openness, satisfying that $\s_k+\tau_k$ 
is $\eta'$-transverse in the compact set $M\cup C\times [1-1/2\alpha,1+1/2\alpha]$. 
Use the $C^1$-openness of the transversality of sections along the distribution $D$ in 
$C$ to assure that $\s_k+\tau_k$ is still $\eta/2$-transverse to {\bf 0} along
$D$ in $C$. This finishes the proof.
\hfill $\Box$

\begin{remark}
Observe that the process followed in the proof is not symmetrical, i.e. we cannot
perturb first the sequence to obtain symplecticity and later on to obtain
contactness, because the perturbations needed to get contactness are not $C^1$
small and so they destroy the achieved simplecticity.

One of the most surprising points of the result is that we cannot
assure $C^0$-closedness between the initial and the perturbed sections. But,
curiously, Donaldson's techniques which are based in this phenomenom continue
applying. For this we will need to control the behaviour of the sections in a certain sense
which will be apparent along the proofs.
\end{remark}

\section{Achieving transversality in local neighborhoods.} \label{borders}
Along this Section we are going to prove Theorem \ref{trans_bor}. We will
assume all the local transversality results developed in \cite{Do99} and \cite{Pr00},
but we need a further refinement to prove the result.

\subsection{Approximately holomorphic models.}
We will use the following Lemma to trivialize the sym-con manifold in the border.
(As before we will enlarge a little the manifold to change the border into a contact
hypersurface). We denote by $C_0$ the subspace of $\R^{2n+2}$ defined as
$$ \{(0,y_0,x_1,y_1, \ldots, x_n, y_n): x_i\in \R, y_j\in \R \}. $$
Moreover we will identify $C_0$ with $\C^n\times \R$ in the natural way.
\begin{lemma} \label{appr_holo}
Given a closed contact manifold $(C, \theta)$ and
a compatible almost-complex structure $J$, construct the $1$-embedding of $C$
into the symplectization $S_D(C)\simeq C\times \R^+$, and denote it by $\hat{C}$. 
Fix a point $x\in \hat{C}$. There exists a uniform constant $c>0$ and a symplectic 
Darboux chart $\varphi: (B_{g}(x,c),\omega) \to (\R^{2n+2}, \omega_0)$ satisfying that:
$\varphi(x)=0$, $\varphi^* J_0(0)=J(x)$, $\varphi_* \hat{D}(x)$ is a complex subspace,
$\varphi^{-1}(C_0)= \hat{C}$ and also
$$ \frac12 g(v,w) \leq \langle (\varphi_*)_yv, (\varphi_*)_yw \rangle \leq
2g(v,w), \forall y\in B_g(x,c),~ v,w\in T_yS_D(C). $$
This implies that $|\nabla^r \varphi|=O(1)$ and $|\nabla^r \varphi^{-1}|=O(1)$,
for r=1,2,3. Also $|\bar{\partial} \varphi(y)|\leq c'd(x,y)$, for a uniform
constant $c'$. \newline
Denote by $\hat{\varphi}$ the restriction of $\varphi$ to $\hat{C}$.
The distribution $\hat{\varphi}_*D$ of $C_0\simeq \C^n\times \R$ can be
equipped with the canonical complex structure $\hat{J}_0$ (obtained by vertical lifting)
and then $|\bar{\partial} \hat{\varphi}(y)|\leq c'd(x,y)$, for all $y\in B_{g}(x,c)
\cap \hat{C}$, where the operator $\bar{\partial}$ is computed respect to $J$ and 
$\hat{J}_0$.
\end{lemma}
{\bf Proof:}
We choose a symplectic Darboux chart at $x$, $\varphi: B_{g_k}(x,c) \to V\subset 
\R^{2n+2}$. The constant $c$ can be chosen in a uniform way because of the compactness 
of $C$. We need to assure also that the standard complex structure $J_0$ in 
$\hat{\varphi}_*D\subset \R^{2n+1}$ and $\hat{\varphi}_*J$ coincide at 
$\hat{\varphi}(x)=0$.
We only have to modify the Arnold's proof of the contact case.
We first use the symplectic Darboux Theorem to obtain Darboux coordinates $\varphi(y)=
(p_0, \ldots, p_n, q_0, \ldots, q_n)$.
 
Following \cite{Ar80} we can assure that the embedding of the contact manifold is
locally given by the equation $p_0=0$. Notice that in general $D_x=\hat{\varphi}_*D(x)
\not = \{ p_0=q_0=0 \}$. But we can choose a standard symplectic basis $(e_1,\ldots,
e_n,f_1, \ldots, f_n)$ in $D_x$. Also we can choose a standard symplectic basis
$(e_0,f_0)$ in $D_x^{\perp}$, assuring that $p_0(e_0)=0$. The orthogonal operation
is made with respect to the symplectic form in the symplectization.
Now, we define the transformation:
\begin{eqnarray}
\eta: \R^{2n+2} & \to & \R^{2n+2} \nonumber \\
\frac{\partial}{\partial p_i} & \to & e_i \nonumber \\
\frac{\partial}{\partial q_i} & \to & f_i. \nonumber
\end{eqnarray}
The map $\eta$ is symplectic and if we compose $\eta \circ \hat{\varphi}$ we obtain
that, in these new Darboux coordinates, denoted again by $(p_0, q_0, \ldots, p_n,q_n)$, 
$C$ is locally defined by the equation
$p_0=0$ and also $D_x$ is complex, in fact $D_x=\{ p_0=q_0=0 \}$. Finally performing
a symplectic transformation in $D_x$ we can assure that $J_{|D}=(J_0)_{|D}$.
Observe that $D_x^{\perp \omega}$ is also complex and then a $Sp(2)$ transformation
there makes that $J(x)=J_0$.

So we have checked that $\varphi^*(J_0)(x)=J(x)$ and $\hat{\varphi}^*(J_0)_{|D}(x)=J_{|D}(x)$
at the point $x$. We cannot assure more because the two complex structures are related
through a, in general non-vanishing, Nijenhuis type tensor at the origin. 
The last inequalities in the statement of the Lemma are assured by the fact
that $\varphi$ is a isometry at $x$ and by the compactness of $C$. Now 
following the discussion in Section 2 of \cite{Do96} it is easy to verify that
the boundings in the antiholomorphic parts are as given. 
\hfill $\Box$

After scaling, the precedent result appears as
\begin{lemma} \label{appr_holok}
Given a closed contact manifold $(C, \theta)$ and
a compatible almost-complex structure $J$. Construct the $1$-embedding of $C$
into the symplectization $S_D(C)\simeq C\times \R^+$, and denote it by $\hat{C}$. 
Fix a point $x\in \hat{C}$. Then there exists a uniform constant $c>0$ and a symplectic 
Darboux chart $\varphi_k: (B_{g_k}(x,c), k\omega) \to (\R^{2n+2}, \omega_0)$ satisfying that:
$\varphi_k(x)=0$, $\varphi_k^* J_0(0)=J(x)$, $(\varphi_k)_* \hat{D}(x)$ is a complex subspace,
$\varphi_k^{-1}(C_0)= \hat{C}$ and also
$$ \frac12 g_k(v,w) \leq \langle ((\varphi_k)_*)_yv, ((\varphi_k)_*)_yw \rangle \leq
2g_k(v,w), \forall y\in B_g(x,c),~ v,w\in T_yS_D(C). $$
This implies that $|\nabla^r \varphi_k|=O(1)$ and $|\nabla^r \varphi_k^{-1}|_{g_k}=O(1)$,
for r=1,2,3. Also $|\bar{\partial} \varphi_k(y)|_{g_k}\leq k^{-1/2}$, for a uniform
constant $c'$. \newline
Denote by $\hat{\varphi}_k$ the restriction of $\varphi_k$ to $\hat{C}$.
Then the distribution $(\varphi_k)_* D$ is a sequence of asymptotically flat contact 
distributions in $\R^{2n+1}$ that are
equipped with the canonical complex structure $\hat{J}_0$ (obtained by vertical lifting)
and then $|\nabla^r \bar{\partial} \hat{\varphi}_k(y)|_{g_k}= O(k^{-1/2})$, for all 
$y\in B_{g_k}(x,c)\cap \hat{C}$ and $r=0,1,2$, where the operator is computed respect 
to $J$ and $\hat{J}_0$.
\end{lemma}
{\bf Proof:}
It follows by composing the map $\varphi$ obtained in Lemma \ref{appr_holo}
with the scaling map $\lambda_k:\C^{n+1} \to \C^{n+1}$ defined as $\lambda_k(z)=
k^{1/2}z$. Then all the boundings are automatic. The only point is to assure
that $|\nabla^r \bar{\partial} \varphi_k(y)|_{g_k}= O(k^{-1/2})$. For $r=0$ it
is a trivial consequence of Lemma \ref{appr_holo}. For $r\geq 1$ follows from 
$|\nabla^r \varphi_k|=O(k^{-(r-1)/2})$. The same occurs with $\bar{\partial}
\hat{\varphi}_k$ and its derivatives.
\hfill $\Box$

\begin{definition}
A sequence of sections $s_k$ of hermitian bundles $E_k$ with
connections has Gaussian decay in 
$C^r$-norm away from the point $x\in M$ if there exists a uniform 
polynomial $P$ and a uniform constant $\lambda>0$ such that for 
all $y\in M$, $|s(y)|$, $|\nabla s(y)|_{g_k}$, $\ldots$, 
$|\nabla^r s(y)|_{g_k}$ are bounded by $P(d_k(x,y))\exp 
(-\lambda d_k(x,y))$. Here $d_k$ is the distance associated 
to the metric $g_k$.
\end{definition}
The following result is used to trivialize bundles in an approximately holomorphic
way. 
\begin{lemma} [\cite{Do96, Au97}] \label{localized}
 Given any point $x\in M$, for $k$ large enough, there exist $(c,c,c)$-bounded 
 sections in $C^r$-norm $\sref$ of $L^{\otimes k}$ over $M$ 
 satisfying the following bounds: 
 $|\sref|>c_s$ at every point of a ball of $g_k$-radius
 $1$ centered at $x$, for some uniform constant $c_s>0$; 
 the sections $\sref$
 have Gaussian decay away from $x$ in $C^r$-norm.
\end{lemma}

\subsection{Some results of approximation theory.}
We give by completeness some basic ideas about the behaviour of the Tchebycheff 
polynomials for interpolating differentiable functions. Finally we prove an
easy, but not standard, result.

In what follows we will study functions $f:[-1,1]\to \C$ and our objective will be to
approximate them by polynomials. We introduce the following
\begin{definition}
The Tchebycheff polynomials $T_n(x)$ are defined inductively as follows
\begin{enumerate}
\item $T_0(x)=1$,
\item $T_1(x)=x$,
\item $T_{n+1}(x)=2xT_n(x)-T_{n-1}(x)$.
\end{enumerate}
\end{definition}
We define the Tchebycheff inner product of two functions $f,g:[-1,1] \to \C$ as
$$ \langle f, g \rangle= \frac2\pi \int_{-1}^1 f \bar{g} \frac{dx}{\sqrt{1-x^2}}. $$
Tchebycheff polynomials satisfy the following simple properties
\begin{lemma} ~ \newline
\begin{enumerate}
\item The system of polynomials $\frac{T_0}{\sqrt{2}}, T_1, T_2, \ldots$ is an orthonormal
system in the space of differentiable functions with respect to the Tchebycheff inner
product.
\item $T_n(x)=\cos(n \arccos x).$
\end{enumerate}
\end{lemma}
{\bf Proof:} It is a simple computation. \hfill $\Box$

Using the precedent result we can compute the orthogonal projection of any given function
to the orthonormal basis $T_0/\sqrt{2}$, $T_1$, etc. So the order $n$ Fourier expansion
of a given function $f$ is
$$ T_nf= \sum_{j=0}^n A_j T_j, $$
where
$$ A_j=\frac{2}{\pi} \int_{-1}^1 f(x) \bar{T_j}(x) \frac{dx}{\sqrt{1-x^2}}. $$
The result we will use is the following technical Lemma, it is nothing but a slight, and
not very precise, adaptation of a classical Jackson's theorem.
\begin{lemma} \label{Tche}
Given $-1<a<b<1$ and given a $C^2$ function $f:[-1,1] \to \C$ which satisfies
that $|f'(x)|<\epsilon$ and $|f''(x)|<\epsilon$ for all $x\in[a,b]$ and $f'(x)=0$, $f''(x)=0$ 
otherwise. Then we have
\begin{enumerate}
\item $|A_j|\leq \frac{4 \epsilon}{\pi j^2}|\arccos b-\arccos a|$
\item $|f-T_nf|_{C^0}\leq \frac{4\epsilon}{\pi n}|\arccos b-\arccos a|$. Therefore 
$T_nf$ converges to $f$ in $C^0$-norm.
\end{enumerate}
\end{lemma}
{\bf Proof:}
The second property follows from the first one by a simple computation summing
the error (and checking that the Fourier expansion converges, which is direct from the 
Weirstrass $M$-test).

To check the first property, we compute it directly. We perform the change of variable
$x=\cos \theta$ and denote $g(\theta)=\bar{f}(cos \theta)$, then we
can write
\begin{equation}
A_j= \frac2\pi \int_0^{\pi} \cos (j\theta) g(\theta)d\theta. \label{coef}
\end{equation}
Integrating by parts,
$$ A_j= \frac{2}{\pi} \int_0^{\pi} \frac1j \sin (j\theta) g'(\theta) d\theta. $$
A new integration by parts leads us to
$$ A_j= \frac{-2}{\pi j^2} \int_0^{\pi}  \cos (j\theta) g''(\theta) d\theta. $$
Now checinkg that $g''(\theta)=\bar{f}''(\cos\theta)\sin^2 \theta- \bar{f}'(\cos
\theta) \cos\theta$ we obtain that $|g''(\theta)|\leq 2\epsilon$ and so
$$ |A_j| \leq \frac{2}{\pi j^2} 2\epsilon \int_{\arccos a}^{\arccos b} 
\cos (j\theta) d\theta \leq \frac{4\epsilon}{\pi j^2} |\arccos b -\arccos a|. $$
So we obtain the required expression. \hfill $\Box$

\subsection{Local result.}
The key point is as usual the local study. We
prove in this Subsection the following
\begin{proposition} \label{trans_cont}
Let $f_k\colon B \times [-1,1] \to \C^m$ be a sequence of functions where $B$ 
is the ball of radius $1$ in $\C^n$ and $B\times [0,1]$ is equipped with 
a sequence of contact forms $\theta(k)$ whose distributions are asymptotically flat. 
Let $0 < \delta < 1/2$ 
be a constant, $\sigma = \delta (\log (\delta^{-1}))^{-p}$, where $p$ is an
integer depending only on the dimensions. 
Assume that $f_k$ satisfies over $B\times  [-1,1]$
the following bounds
\begin{eqnarray*}
|f_k| \leq 1, ~~~~ |\bar{\partial}_0 f_k| \leq \sigma, ~~~~|\nabla
\bar{\partial}_0 f_k| \leq \sigma, \\
|\partial f_k /\partial s| < 1, ~~ |\partial \nabla f_k /\partial
s| < 1, 
\end{eqnarray*}
for $k$ large enough, where $\bar{\partial}_0$ is the $(0,1)$ operator defined in 
$D(k) = \ker\theta(k)$ by the complex structure $J_0$ and $s$ is the real coordinate. 
Then for $k$ large enough
there exists a holomorphic polynomial $t_k: \C \to \C^m$  such that 
$|t_k|< \delta$ on the set $[-2k^{1/6}, 2k^{1/6}]\times \{0 \}\subset \C$
and such that the function $s_k(z_1,\ldots, z_{n+1})=f_k(z_1, \ldots, z_{n+1})-t_k(z_{n+1})$ is 
$\sigma$-transverse along the distribution $D(k)$ to zero 
on $B(0,1/2)\times [-1,1] \subset \C^n \times \C=\C^{n+1}$ for $k$ large enough. 
Moreover, the modulus of $t_k$ and of its first and second derivatives can be bounded above 
by a fixed real polynomial $b_{\delta}$ depending only on $\delta$.
\end{proposition}
This Proposition is a consequence of the local transversality results of the contact
category which are stated in all generality in \cite{Pr00} as
\begin{proposition}[Proposition 4.4 in \cite{Pr00}] \label{trans_ant}
Let $f_k\colon B \times [0,1] \to \C^m$ be a sequence of functions where $B$ 
is the ball of radius $1$ in $\C^n$ and $B\times [0,1]$ is equipped with 
a sequence of contact forms $\theta(k)$ whose distributions are asymptotically flat. 
Let $0 < \delta_0 < 1/2$ 
be a constant and let $\sigma = \delta_0 (\log (\delta_0^{-1}))^{-p}$, where $p$ is a
integer depending only on the dimensions. Assume that $f_k$ satisfies over $B\times  [0,1]$
the following bounds
$$ |f_k| \leq 1, ~~~~ |\bar{\partial}_0 f_k| \leq \sigma, ~~~~|\nabla
\bar{\partial}_0 f_k| \leq \sigma,$$
for $k$ large enough, where $\bar{\partial}_0$ is the $(0,1)$ operator defined in $D(k) = \ker
\theta(k)$ by vertical projection of the standard complex structure $J_0$.  Then 
for $k$ large enough
there exists a smooth curve $w_k\colon [0,1] \to \C^m$  such that $|w_k|< \delta_0$
and the function $f_k - w_k$ is $\sigma$-transverse to zero on $B(0,1/2)\times [0,1]$. 
Moreover, if $|\partial f_k /\partial s| < 1$ and $|\partial \nabla f_k /\partial
s| < 1$, we can choose $w_k$ such that $|d^iw_k /ds^i | < \Phi(\delta_0)$, $(i=1,2)$; 
$d^jw_k /ds^j (0) = 0$ and 
$d^jw_k / ds^j (1) = 0$, for all $j \in \N$, where $c$ is a uniform  constant and $\Phi:\R^+
\to \R^+$ is a function depending only on the dimensions.
\end{proposition}
{\bf Proof of Proposition \ref{trans_cont}:}
Our hypothesis coincide with the ones in Proposition \ref{trans_ant}. 
We choose $\delta_0=\delta/2$. So we obtain a function $w_k:[-1,1] \to \C^m$, such that 
$f_k-w_k$ is $\sigma$-transverse and satisfying also that $|w_k|\leq \delta/2$.
The idea is to approximate $w_k$ by a complex polynomial. First we extend $w_k$ to
the whole real line as
$$ \hat{w}_k(t)= \left\{ \begin{array}{ll}
w_k(-1) & \mbox{if} ~ x\leq -1. \\
w_k(t)& \mbox{if} ~ -1\leq x\leq 1. \\
w_k(1) & \mbox{if} ~ x\geq 1.
\end{array}
\right. $$
Now we scale the real coordinate constructing a new function $h_k(x)=\hat{w}_k(2k^{1/6}x)$. 
Obviously we have the following boundings $|h_k(x)|\leq \delta$, 
$|\frac{d^j h_k}{ds^j}|\leq 2k^{1/6}\Phi(\delta)$ for $j=1,2$.
Moreover $|\frac{d^j h_k}{ds^j}|=0$ if $x\in [-\frac{1}{2k^{1/6}}, \frac{1}{2k^{1/6}}]$.

Decompose $h_k=(h_k^1, \ldots, h_k^m)$. Then each of the components $h_k^j$ is
in the hypothesis of Lemma \ref{Tche}, when restricted to the segment $[-1,1]$. So we have 
that the associated Tchebycheff polynomial of degree $d$ satisfies 
$$|h_k^j-T_dh_k^j|_{C^0}\leq \frac{8\Phi(\delta)k^{1/6}}{\pi d}
|\arccos(\frac{1}{2k^{1/6}})- \arccos(-\frac{1}{2k^{1/6}})|.$$
We substitute $\pi/2-|x| \leq |\arccos (x)| \leq \pi/2+|x|$.
Summing up all the components we find
$$|h_k-T_dh_k|_{C^0}\leq \frac{8m\Phi(\delta)}{\pi d}, $$
Then increasing enough $d$ we can assure that 
\begin{eqnarray}
|h_k-T_dh_k|\leq \frac{\sigma}{2c_u}, \label{sisi} \\
|h_k-T_dh_k|\leq \frac{\delta}{2}, \label{sisisi}
\end{eqnarray}
In fact, we need $d=O(\max \{ \sigma^{-1}\Phi(\delta), \delta^{-1} \})$, where $c_u$ is 
the uniform constant
of $C^1$-openness for the property of being transverse to {\bf 0} along the distribution
$D(k)$. Define $t_k(z)=T_dh_k(\frac{z}{2k^{1/6}})$.
So, we claim that, imposing (\ref{sisi}) and (\ref{sisisi}), $f_k-t_k$ is 
$\sigma/2$-transverse
to {\bf 0} along $D$ in $B\times [-1,1]$ as we wanted, and also that $|t_k(z)|\leq \delta$
for all $x\in [-2k^{1/6},2k^{1/6}]\times \{ 0 \}$. 
To prove it, we extend to $\C^n\times \R$ the functions $w_k$ and
$t_k$ as $w_k(z_1,\ldots,z_n,s)=w_k(s)$ and $t_k(z_1,\ldots,z_n,s)=t_k(s)$. Now it is
easy to check that in $B\times [-1,1]$ the function $w_k-t_k$ has, for k large, 
mixed $C^2$-boundings $(\frac{\sigma}{2c_u}, c_R)$, where $c_R$ is a constant depending 
only on $\delta$. Therefore, recalling that $f_k-w_k$
is $\sigma$-transverse along $D$, we obtain that $f_k-t_k$ is $\sigma/2$-transverse
to $D$ in $B\times [-1,1]$.

To finish we need to bound above the modulus of $t_k$, or equivalently $T_dh_k$, by a 
fixed polynomial. For this we need only to recall the first property of
Lemma \ref{Tche} which translates in our case
$$|\hat{A}_j^l|\leq \frac{4 \Phi(\delta)}{\pi j^2},$$
where $\hat{A}_j^l$ is the $A_j$ component of the polynomial $T_dh_k^l$ once the rescaling
$2k^{1/6}$ is introduced.
It implies that the coefficient $\hat{A}_j$ is bounded above by a function of $\delta$.
So, for a fixed $\delta$ the degree $d$ is constant and the coefficients of the 
Tchebycheff aproximation are bounded above by a constant. Then it is obvious that
there exists a fixed real polynomial bounding above the modulus of $t_k$ and of its 
derivatives. This finishes the proof.
\hfill $\Box$

The following result has a more geometrical appearance.
\begin{proposition} \label{local_sol}
Let $C$ be a cooriented contact manifold and let $s_k$ be a sequence of sections with 
global $C^3$-boundings $(c,c_D,c_R)$ of the bundles $E\otimes L^{\otimes k}$ over the 
symplectization $S_D(C)$.
Then given a point $x$ in the $1$-embedding $\hat{C}$ and $\delta>0$ there exists a 
sequence of sections 
$\tau_{k,x}$ of $E\otimes L^{\otimes k}$ and $\sigma = \delta(\log (\delta^{-1}))^{-p}$ 
(for some integer $p>0$) satisfying that:
\begin{enumerate}
\item $\tau_{k,x}$ has global $C^3$-boundings for $k$ large (depending on $\delta$)
\begin{eqnarray*}
& & (c_u c_R P_{\delta}(d_k(x,y))
\exp(-\lambda d_k(x,y)^2), c_u c_D  \delta Q(d_k(x,y))\exp(-\lambda d_k(x,y)^2), \\
& & c_uc_R P_{\delta}(d_k(x,y))\exp(-\lambda d_k(x,y)^2))
\end{eqnarray*}
at any point $y$,
\item $(s_k+\tau_{k,x})_C$ is $\sigma$-transverse to {\bf 0} along $D$ in 
$B_{g_k}(x,\hat{c})\bigcap \hat{C}$
\end{enumerate}
for $k$ large enough, where $\lambda$ and $p$ are constants depending only on the dimensions, 
$P_{\delta}$ is a uniform polynomial (depending on $\delta$), $Q$ is a uniform
polynomial (not depending on $\delta$), $\hat{c}$ and $c_u$ are 
uniform constants.
\end{proposition}
{\bf Proof:}
We choose the trivializations $\varphi$ and $\varphi_k$ defined by Lemmas
\ref{appr_holo} and \ref{appr_holok}. Also we fix a section $s_{k,x}^{ref}$ as defined
in Lemma \ref{localized}. Fix a unitary basis $\{e_1(x),\ldots, e_r(x)\}$ in $E_x$ and
extend it by parallel transport along radial directions to a frame
$\{e_1, \ldots, e_r\}$ in a neighborhood of $x$. It is easy to check that
$|\nabla^r e_i|_{g_k}=O(k^{-r/2})$ and so the sequence of sections $a_k^j=e_j$
has $c$ bounding in $C^r$-norm, for a uniform $c>0$. Now we define the
frame:
$$ \s_j=e_j \otimes s_{k,x}^{ref}, $$
which is bounded in $C^r$-norm by construction. Moreover $|\s_j|>c_s$ for any
$y\in B_{g_k}(x,1)$. Finally choosing a sufficiently small uniform $\hat{c}$,
we have that $\s_1, \ldots, \s_r$ is approximately unitary for any 
$y\in B_{g_k}(x,\hat{c})$.

Now we construct an application $\tilde{f}_k:B_{g_k}(x,\hat{c}) \to \C^r$ imposing 
the condition
$$ s_k(y)= \tilde{f}_k^1(y)\cdot \s_1(y) +\ldots \tilde{f}_k^r(y) \cdot \s_r(y). $$
Using that $\s=(\s_1, \ldots, \s_r)$ is approximately unitary, namely, interpreted in 
each fiber
as a linear application $\s:\C^r \to E_y$, $\s$ has an inverse with uniformly bounded norm,
we find
\begin{equation}
|\tilde{f}_k|\leq c_u, ~~ |\nabla^r \tilde{f}_k|\leq c_u, ~~ | \nabla^{r-1} 
\bar{\partial} \tilde{f}_k|\leq c_uk^{-1/2}, \label{cotass}
\end{equation}
for $r=1,2,3$. 
Finally we use the chart $\varphi_k$ to define an application $f_k= \tilde{f}_k 
\circ \varphi_k^{-1}$. Scaling the chart by an appropiate uniform constant
we can assure that $\varphi_k(B_{g_k}(x,\hat{c}/8))\subset B(0,1/2) \subset B(0,2)
\subset \varphi_k(B_{g_k}(x, \hat{c}))$. This is possible, perhaps after shrinking
uniformly $\hat{c}$, because of the approximately isometry property of Lemma
\ref{appr_holok}. From (\ref{cotass}) and the boundings of Lemma \ref{appr_holok} we obtain
\begin{equation}
|f_k|\leq c_u, ~~ |\nabla^r f_k|\leq c_u, ~~ | \nabla^{r-1} 
\bar{\partial} f_k|\leq c_uk^{-1/2}, \label{cotas2}
\end{equation}
Without loss of generality we suppose that $f_k$ satisfies the boundings required
in Proposition \ref{trans_cont} (in fact, we only have to multiply it by a non-zero uniform
constant to assure this). Then the precedent Proposition applies, once $k$ is
large enough, and we obtain a polynomial $t_k$ satisfying the conditions of Proposition
\ref{trans_cont}. We extend the definition of $t_k$ to $\C^{n+1}$ as
$$ t_k(z_1, \ldots,z_{n+1})=t_k(z_{n+1}). $$
Now we define $\tilde{t}_k= t_k \circ \varphi_k$. Recall that this is defined in a ball
of $g_k$-radius $O(k^{1/2})$, which is the domain of $\varphi_k$ (obviously, it has the
same domain than $\varphi$). Then, taking into
account that $s_{k,x}^{ref}$ has support in a ball of $g_k$-radius $O(k^{1/6})$, 
define
$$ \tau_{k,x}= \tilde{t_k}^1 \s_1+ \ldots \tilde{t_k}^r \s_r. $$
Recall that $s_k+\tau_{k,x}$ and $f_k+t_k$ are related through uniform scaling
constants and through the approximately unitary basis $\s$. So, by construction,
the property $2$ of the statement is satisfied, except by a uniform multiplying factor
which can be eliminated by increasing uniformly the integer $p$.

Recall that we can bound $t_k$ by a fixed polynomial $b_{\delta}$.
Define $\hat{P}_{\delta}(r)= \max_{|z|=r} \{ b_{\delta}(r) \}$.  We can find
a fixed polynomial $P_{\delta}(r)$ satisfying that $P_{\delta}(r)\geq \hat{P}_{\delta}(r)$.
Then we 
easily conclude that $\tau_k$ has global boundings 
\begin{eqnarray*}
& & (c_u c P_{\delta}(d_k(x,y))
\exp(-\lambda d_k(x,y)^2), c_u c_D  \delta Q(d_k(x,y))\exp(-\lambda d_k(x,y)^2), \\
& & c_uc_R P_{\delta}(d_k(x,y))\exp(-\lambda d_k(x,y)^2)),
\end{eqnarray*} 
for any $y\in \hat{C}$. 
The first and third boundings are trivial. For the second one we proceed
as follows. The bounding of $|\tau_{k,x}(y)|$ follows form the condition
that $|t_k(y)|<\delta$ for all the points of the set $\C^n \times [-2k^{1/6},
2k^{1/6}] \times \{0 \}$, which implies that $|\hat{t}_k(y)|<c_u\delta$
at any point $y\in \hat{C}\cap B_{g_k}(x, k^{1/6})$. For bounding the derivatives,
we denote $D_0$ the pull-back through the map $\hat{\varphi}_k$ of the distribution 
$\C^n \times \{ 0 \}$ . We easily bound
\begin{equation}
\angle_M(D(y), D_0(y)) \leq ck^{-1/2}d_k(x,y), \label{cerquita}
\end{equation}
where $c>0$ is certain uniform constant. By construction,
\begin{equation}
\nabla_{D_0} \hat{t}_k= 0. \label{inmovil}
\end{equation}
And so using (\ref{cerquita}), (\ref{inmovil}) and the bounding polynomial
$P_{\delta}$ we find
$$ |\nabla_{D} \hat{t}_k(y)|= ck^{-1/2}d_k(x,y)P_{\delta}(d_k(x,y)). $$
We change the polynomial $P_{\delta}(t)$ by $t\cdot P_{\delta}(t)$ and so
$$ \nabla_{D} \hat{t}_k(y)= ck^{-1/2}P_{\delta}(d_k(x,y)). $$
Therefore
\begin{eqnarray*}
|\nabla_D (\hat{t}_k\cdot \sigma)| & = & |\nabla_D \hat{t}_k \cdot \sigma +
\hat{t}_k \cdot \nabla_D \sigma| = \\
& \leq & ck^{-1/2}P_{\delta}(d_k(x,y))\exp(-\lambda d_k(x,y)^2)+ \\
& + & \delta Q(d_k(x,y))\exp(-\lambda d_k(x,y)^2), 
\end{eqnarray*}
which, for $k$ large enough, satisfies the required bounding because the first term
is arbitrarily small and the polynomial $Q$ does not depend on $\delta$ as required.
The boundings on $|\nabla^r \tau_{k,x}|$ are obtained in the same way.
\hfill $\Box$

\subsection{Globalization process.}
As in \cite{Do96, IMP99} the final point will be to construct a global perturbation
of the sequence of sections from a sequence of localized perturbations added
in a suitable way. Along this Subsection we adapt Donaldson's framework to our case.
This development is given by

\noindent {\bf Proof of Theorem \ref{trans_bor}:}
Donaldson's globalization argument works with some slight
variations. Choose a finite set of points $S$ satisfying the following conditions:
\begin{enumerate}
\item $\cup_{x\in S} B_{g_k}(x,\hat{c}) \supset \hat{C}$.
\item There exist a partition $S=\cup_{j\in J} S_j$ verifying that $d_{g_k}(x,y)>
N$ if $x,y\in S_j$. $N$ will be fixed along the proof.
\item The cardinal of $J$ is $O(N^{2n+1})$.
\end{enumerate}
Recall that the starting sequence of sections has global $C^3$-boundings $(c,c_D,
c_R)$. We proceed by steps perturbing at each $S_j$ at a time.
Let us find a perturbation centred on each of the points of $S_1$ to achieve
trasnversality at a neighborhood of $S_1$. Fix $x\in S_1$, use Proposition
\ref{local_sol} with certain $\delta=\delta_1>0$ to be chosen. We find out a sequence
of perturbations $\tau_{k,x}$ with global boundings 
\begin{eqnarray*}
& & (c_u c_R P_{\delta}(d_k(x,y))
\exp(-\lambda d_k(x,y)^2), c_u c_D  \delta Q(d_k(x,y))\exp(-\lambda d_k(x,y)^2), \\
& & c_uc_R P_{\delta}(d_k(x,y))\exp(-\lambda d_k(x,y)^2)).
\end{eqnarray*}
We take $\delta_1$
to assure that the second bounding is uniformly less that $\epsilon/2$ (Recall
that $\epsilon$ is the maximum first mixed $C^3$-bounding admitted). In fact,
we can choose $\delta_1\leq c_p \epsilon$, for a certain uniform $c_p>0$. This
is possible  since $Q$ does not depend on $\delta$! 
We have now a perturbation centred on each of the points of $S_1$ which solves
the problem in the balls $B_{g_k}(x,\hat{c})$. But the perturbations are not
independent. This is the moment when the integer $N$ comes into play.
Analyze a fixed $x\in S_1$. We can compute
the maximum first mixed $C^3$-bounding of the perturbations (the boundings in the distribution
directions) of the rest of the points of $S_1$ in the ball
$B_{g_k}(x,c)$. This ``bad'' perturbation is bounded by $c_u\delta\exp(-\lambda N^2)$.
Again, it is very important to assure that $Q$ does not depend on $\delta$
to find $c_u$ independent of $\delta$, otherwise the globalization process does not hold.
To avoid the destruction of the achieved transversality a sufficient condition is so
\begin{equation}
 c_u\delta_1\exp(-\lambda N^2) \leq \delta_1 (\log (\delta_1^{-1}))^{-p}, \label{aguanta}
\end{equation}
for a uniform constant $c_u$ not depending on $\delta$. In this first
stage we may choose N to satisfy (\ref{aguanta}).
So adding all the perturbations we find a sequence of sections $\tau_k^1$ which
added to $s_k$ achieve $\sigma_1$-transversality in $\bigcup_{x\in S_1}
B_{g_k}(x,\hat{c})\cap \hat{C}$. Moreover we find that the sequence $\tau_k^1$ has global
$C^3$-boundings $(c_1', \epsilon/2, c_2')$. We only know that $c_1'$ and $c_2'$
depend on $\delta$, but the important point is that ``they exist''. 
Now in the second stage we choose $\delta_2$ to assure that the final sequence of
perturbations $\tau_k^2$ has boundings $(c_2', \min \{ \epsilon/4, \sigma_1/{2c_u} \}, c_2')$.
The second bounding is imposed to guarantee that the sequence has controled boundings in 
the distribution $D$ directions and also that do not destroy the achieved transversality
in the $\hat{c}$-neighborhood of $S_1$ ($c_u$ is the constant of $C^1$-openness
of the transversality to {\bf 0} along $D$).

Repeating the process we find $\tau_k=\sum_{j=1}^{q} \tau_k^j$ that
has global $C^3$-boundings $(c', \epsilon, c_R')$, which are independent of $k$ because
$q$ is independent. Moreover $s_k+\tau_k$ is $\sigma$-transverse to {\bf 0}
along $D$ all over $\hat{C}$. Again, the constant $\sigma>0$ is uniform because
the number of steps is independent of $k$.

Only one important question has to be checked. The expression (\ref{aguanta})
must hold in all the steps of the process. Namely we must assure
$$ c_u\delta_j\exp(-\lambda N^2) \leq \delta_j (\log (\delta_j^{-1})^{-p}. $$
But, the asymptotic analysis of the expresion $(\log (\delta_j^{-1})^{-p}$
provides this condition if we choose $N$ large enough (for a proof of this fact see
Section 2 in \cite{Do96}).
\hfill $\Box$

\section{Topological considerations.} \label{topol}
In this Section we characterize the topological properties of the constructed submanifolds.

\subsection{Relative Lefschetz hyperplane theorem.}
We prove now the second part of Theorem \ref{main_thm}. The started point is a 
sym-con manifold $(M,\omega, C,\theta)$ where we have found a sequence of sym-con
submanifolds $(W_k, \omega, C_k, \theta)$ obtained as zero sets of a sequence
of sections $s_k$ of the bundles $E\ox L^{\ox k}$ with global boundings which are transverse 
to {\bf 0} in the symplectic manifold and in the contact border.

We take as a tool the functions $f_k(p)=\log |s_k(p)|^2$. Then we follow the
Proof of Proposition 2 in Section 5.1 of \cite{Au97} to conclude that
the critical points of $f_k$ in $M$ have at least index $n-r+2$, for $k$
large enough. In the same way following \cite{IMP99} we conclude that the
critical points of $f_k$ in the border $C$ have at least index $n-r+1$ (again, for
$k$ sufficiently large).
So, being the border a closed manifold, this proves that the inclusion
$i:C_k \to C$ induces isomorphism in homology groups $H_j$ (resp. homotopy groups)
for $j\leq n-r-1$ and surjection for $j=n-r$. This is the content of the Lefschetz
theorem in the contact case.

We are going to define the double
copy $M^d$ of the manifold $M$ as the topological connected sum $M\cup_C M$. We can 
arrange this topological operation, for each $k$, to assure the smoothness of the 
submanifold $W_k^d=W_k\cup_{C_k} W_k$.

Now we perturb $f_k$ into a new function $\hat{f}_k$ to assure that the natural extension to 
the double copy is smooth. For this we only need to assure 
that $\frac{df_k}{n}(c)=0$ for any $c\in C$ where $n$ is the normal direction to $C$ respect to 
the metric $h_k$. Use that in a small neighborhood $V$ of $C$ we can trivialize $M$ as 
$C\times [0,\epsilon)$ assuring also that $n=\frac{\partial}{\partial s}$ being $s$ the 
real coordinate. Therefore we perturb $f_k$ in this small meighborhood as
$$ \hat{f}_k(c,s)=f_k(c, \beta(s)), $$
where $\beta:[0,\epsilon] \to [0,\epsilon]$ is a smooth function satisfying
\begin{enumerate}
\item $\beta(0)=0$ and $\beta(\epsilon)=\epsilon$.
\item $\beta'(x)>0$ for all $x\in (0,\epsilon)$.
\item $\frac{d^r \beta(0)}{ds^r}=0$, for all $r\in \N^*$.
\item $\beta'(\epsilon)=1$ and $\frac{d^r \beta(\epsilon)}{ds^r}=0$ for $r=2,3, \ldots$
\end{enumerate}
It is easy to check that $\hat{f}_k$ extends to a smooth function, again denoted,
$\hat{f}_k$ in $M^d$. Moreover the critical points of $f_k$ and $\hat{f}_k$ coincide
in $M$ because $\beta$ only performs a diffeomorphism
outside the border. The indices do not change. In $C$ (interpreted as a submanifold in 
$M^d$) we obtain that the critical points of $f_k$ are now critical points of
$\hat{f}_k$ and the index of these critical points is at least $n-r+1$.

Summarizing, the manifold $W_k^d$ is a smooth submanifold of $M^d$.
The function $\hat{f}_k$ has critical points of index at least $n-r+1$ for $k$ large enough.
This implies, by standard Morse theory, that the inclusion
$$ i_d: W_k^d \to M^d $$
induces isomorphisms in homology and homotopy groups for dimension less than or equal to
$n-r$ and surjection for $n-r+1$. 

We denote the first and second copies of $M$ in $M^d$ by $M^1$ and $M^2$
respectively. The same for $W_k$ with copies $W_k^1$ and $W_k^2$. The natural
diffeomorphism defined in $M^d$ interchanging the copies is denoted
as $e:M^d \to M^d$, namely $e(M^1)=M^2$, $e(M_2)=M_1$ and $e(C)=C$.
 
Our objective is to prove that the natural morphism
$$i_j: H_j(\bar{W}_k, C_k) \to H_j(\bar{M},C) $$
is actually an isomorphism when $j\leq n-r$ (and an epimorphism for $j=n-r+1$). 
There are several ways to
prove this result we choose a constructive one which is a little longer than others
because it clarifies a bit the topological ideas involved in the proof.

First let us prove that $i_j$ is epimorphism in the required cases.
We choose $\alpha^1\in H_j(\bar{M},
C)$. Identify $M\simeq M^1$, then we define $\alpha^2=-e_* \alpha_1$. Construct
$\alpha^d=\alpha_1+\alpha_2$ which is an element of $H_j(M^d)$. If $j\leq n-r+1$,
then there exists $\gamma^d\in H_j(W_k^d)$ such that $\alpha^d-\gamma^d=\partial \epsilon$,
for some $\epsilon\in H_{j+1}(M^d)$.

After a small isotopic perturbation, we can suppose that all the elemental chains defining 
$\gamma^d$ and $\epsilon$ are trasnverse to $C$. We claim that we can find an element
homologous to $\gamma^d$ of the form 
\begin{equation}
\gamma^1+f+\gamma^2-f, \label{deco}
\end{equation}
where $\gamma^i\in H_j(\bar{W}_k^i, C_k)$, $i=1,2$ and $f$
is a chain in $C$. For this recall that given an
elemental $1$-chain $c:[0,1]\to M^d$ transverse to $C$ we can define a $1$-chain $c^1$ in
$M^d$ as follows. The chain intersects $c$ in points $a_1,b_1,a_2,\ldots$ (suppose
that $c(0)\in M^1$). Then we define $c^1=c[0,a_1]+c[b_1,a_2]+\cdots$. The same hold
for any $j$-chain transverse to $C$ using an adequate triangulation. In fact, the morphism
$c\to c^1$ is the explicit way of the composition
$$ H_j(W_k^d) \to H_j(W_k^d, \bar{W}_k^2) \to H_j(\bar{W}_k^1, C_k), $$
where the first morphism is the restriction and the second one is generated by excision
of $W_k^2$. Then by the construction we find that $(\alpha^d)^1=\alpha^1$ and
$\alpha^1-\gamma^1-f=\partial \epsilon^1$, for certain chain $f$ in $C$.
This implies (\ref{deco}) and $\alpha_1$ and $\gamma^1$ are homologous relative to
the border $C$. So the morphism $i_j$ is surjective in the expected cases.

Now we study the injectivity. Choose $\gamma^1\in H_j(\bar{W}_k, C_k)$ and suppose that
there exists a $(j+1)$-chain $\epsilon^1$ in $\bar{M}$ such that $\partial \epsilon^1 = \gamma^1-
c$, for some chain $c$ of $C_k$. The question is whether we are able to
find a chain $\epsilon^1$ in $\bar{W}_k$. The argument is analogous to the precedent one.
We construct $\gamma^2= -e_* \gamma^1$ and therefore $\gamma^d=\gamma^1+\gamma^2$
is an element of $H_j(W_k^d)$. In the same way we construct $\epsilon^d$ satisfying
$\partial \epsilon^d= \gamma^d$. If we assume that $j\leq n-r$ then the Lefschetz
hyperplane theorem in $M^d$ assures that there exists a $(j+1)$-chain $\rho$ in $W_k^d$
satisfying $\partial \rho= \gamma^d$. Now we construct $\rho^1$, which is in $\bar{W}_k$, 
as in the precedent case. Finally, we obtain $\partial \rho^1 =\gamma^1 -c$, for 
some chain $c$ in $C_k$.
\hfill $\Box$

\subsection{Homology and Chern numbers of the submanifolds.}
To finish we state the following straightforward result.
\begin{proposition}
Given any sequence of sections $s_k$ with global $C^3$-boundings of bundles $E\otimes
L^{\ox k}$ which are transverse to zero, then the Chern classes of the symplectic
zero sets $Z(s_k)$ are given by
$$ c_l(TZ(s_k))=(-1)^l \left( \begin{array}{c} r+l-1 \\ l \end{array} \right)
(k \left[ \frac{\omega}{2\pi}\right] )^l+O(k^{l-1}). $$
\end{proposition}
{\bf Proof:}
Denote $Z(s_k)=W_k$. The formula follows directly from the relation
$$ i^* c(TX)= i^* c(E\otimes L^{\otimes k}) \cdot c(TW_k). $$
\hfill $\Box$


\begin{thebibliography}{xxxxx}

\bibitem[Ar80]{Ar80}  V. Arnold.  Mathematical Methods of Classical Mechanics.  
Springer-Verlag (1980).
 
\bibitem[Au97]{Au97}  D. Auroux.  {\it Asymptotically holomorphic families of
symplectic submanifolds}.  Geom. Funct. Anal., {\bf 7}, 971-995 (1997).

\bibitem[Au99]{Au99}  D. Auroux.  {\it Th\'eor\`emes de structure des vari\'et\'es 
symplectiques compactes via des techniques presque complexes}.  Ph. D. Thesis.  (1999).

\bibitem[Au99b]{Au99b} D. Auroux, L.  Katzarkov. {\it Branched coverings of $\CP^2$ and
associated invariants of symplectic $4$-manifolds.} Preprint. \'Ecole
Polytechnique. Paris. (1999).

\bibitem[Ch66]{Ch66} E. W. Cheney. {\it Introduction to Approximation Theory.}
McGraw-Hill Book Company. New York. (1966).

\bibitem[Do96]{Do96}  S. K. Donaldson. {\it Symplectic submanifolds and
almost-complex geometry}.  J. Diff. Geom., {\bf 44}, 666-705 (1996).


\bibitem[Do99]{Do99}  S. K. Donaldson. {\it Lefschetz pencils on symplectic manifolds}. 
Preprint (1999).



\bibitem[El98]{El98}  Y. Eliashberg.  {\it ICM 98, Berlin}  (1998).



\bibitem[IMP99]{IMP99} A. Ibort, D. Mart\'{\i}nez, F. Presas. {\it On the construction of
contact submanifolds with prescribed topology}. Preprint. Universidad Carlos III de Madrid.
(1999).






\bibitem[MPS99]{MPS99} V. Mu\~noz, F. Presas, I.Sols. {\it Almost holomorphic embeddings
in grassmanians with applications to singular symplectic submanifolds.} Preprint. Universidad
Complutense de Madrid. (1999)



\bibitem[Pr00]{Pr00} F. Presas. {\it Lefschetz type pencils on contact manifolds.}
Preprint, Universidad Complutense de Madrid. (2000).




\end{thebibliography}
\end{document}